
\input amstex
\documentstyle{amsppt}

\NoRunningHeads
\magnification=1100
\baselineskip=29pt
\parskip 9pt
\pagewidth{5.2in}
\pageheight{7.2in}

\TagsOnLeft
\CenteredTagsOnSplits

\def\alp{\alpha}
\def\and{\text{and}}

\def\bl{\bigl(}
\def\br{\bigr)}
\def\Bl{\Bigl(}
\def\Br{\Bigr)}

\def\cc{\Cal C}
\def\ch{\text{ch}}
\def\co{\tilde c_1}
\def\ct{\tilde c_2}
\def\ctf{\tilde c_2(\f)}
\def\CC{\bold C}
\def\coker{\operatorname{coker}}
\def\codim{\operatorname{codim}}

\def\det{\operatorname{det}}
\def\dual{^{\vee}}

\def\e{\Cal E}

\def\f{\Cal F}
\def\ff{\pab\sta\f\tor\pac\sta\f}

\def\g{\Cal G}

\def\half{{1\over 2}}

\def\lra{\longrightarrow}
\def\lreg{_{\text{reg}}}
\def\lowq{_{\QQ}}
\def\lowo{_{[1]}}
\def\lowf{_{[5]}}

\def\md{\frak M(I,d)}
\def\mbd{\overline{\frak M}(I,d)}
\def\mthd{{\frak M}_{\tilh}(I,d)}
\def\mhdb{\overline{\frak M}_H(I,d)}
\let\mdb=\mbd
\def\mhd{\frak M_H(I,d)}
\def\mh{\!:\!}
\def\m{\frak M}
\def\mapright#1{\,\smash{\mathop{\lra}\limits^{#1}}\,}

\def\nqp{N_{\QQ}^+}

\def\OO{\Cal O}
\def\otime{\!\otimes\!}
\def\oo{\Cal O(1)}

\def\PP{\bold{P}}
\def\pix{\pi_X}
\def\picff{\pbc\sta C\tor\pab\sta\f\tor\pac\sta\f}
\def\picxs{\pic\bl\xxx\br^{\sig}}
\def\pri{^{\prime}}
\def\pa{p_1}
\def\pb{p_2}
\def\pc{p_3}
\def\pab{p_{12}}
\def\pbc{p_{23}}
\def\pac{p_{13}}
\def\pe{\bold P(\e)}

\def\pia{\pi_1}
\def\pib{\pi_2}
\def\pic{\text{Pic}}
\def\picmd{\pic\bl \md\br}
\def\picmdb{\pic\bl \mdb\br}
\def\plgg{\pbc\sta\psi(L)\tor\pab\sta\g\tor\pac\sta\g}
\def\pr{\operatorname{pr}}

\def\qu{\bold{Quot}}
\def\quot{/\!/}
\def\qxx{q_{XX}}

\def\ss{^{ss}}
\def\sta{^{\ast}}
\let\sig=\sigma
\def\stab{\text{stab}}
\let\sub=\subseteq
\def\subb{\Subset}

\def\tilx{\widetilde{X}}
\def\tilh{\tilde H}
\def\time{\!\times\!}
\def\td{\operatorname{td}}
\def\timezq{\otimes_{\ZZ}\QQ}
\def\tdxx{\td(X\time X)}
\def\tdx{\td(X)}
\def\time{\!\times\!}
\def\tor{\boxtimes}
\def\tilk{\widetilde{K}}
\def\tz{\!\otimes_{\ZZ}\ZZ[{1\over 12}]}
\def\tpa{\tilde p_1}

\def\tpab{\tilde p_{12}}
\def\tpbc{\tilde p_{23}}
\def\tpac{\tilde p_{13}}

\def\umo{^{-1}}

\let\vphi=\varPhi
\def\vphib{\overline{\vphi}}

\def\xxx{X\times X}
\def\xx{\chi}

\def\ZZ{{\Bbb Z}}
\def\QQ{{\Bbb Q}\,}
\def\RR{{\Bbb R}}
\def\CC{{\Bbb C}}
\let\pro=\proclaim
\let\endpro=\endproclaim
\def\proof{\noindent{\bf Proof}.\ }

\topmatter
\title
Picard groups of the moduli spaces of vector bundles over
algebraic surfaces
\endtitle

\author
Jun Li
\endauthor
\thanks
This work was partially supported
by NSF grant DMS-9307892 and A. Sloan fellowship
\endthanks

\affil
Mathematics Department\\
University of California, Los Angeles
\endaffil

\email
jli@math.ucla.edu
\endemail

\endtopmatter
\document

The purpose of this note is to determine the Picard group of the
moduli space of vector bundles over an arbitrary algebraic surface.
Since Donaldson's pioneer work on using moduli of vector bundles to
define smooth invariants of an algebraic surface, there has been a surge
of interest in understanding the geometry of
this moduli space. Among other things, the study of line bundles
on this moduli space plays a major role
in this area. One important question remain open until now is to
determine the Picard group of
this moduli space. This is known in some special cases, for instance
for projective plane [St], ruled surfaces [Qi2, Yo] and K3 surfaces [GH].
In this note, we will settle this question by
providing a general construction of line bundles that will
include virtually all line bundles on this moduli space, when the second
Chern class of the sheaves parameterized by this moduli space is
sufficiently large.
The construction is again based on Knudsen and Mumford's
recipe of determinant line bundle construction. The new input from this note
is that instead of using complexes on surface $X$\
we will use complexes on $\xxx$, which yields
all previously known line bundles as well as new ones.
After this construction,
we will use our knowledge of the first two Betti numbers of
this moduli space to argue that this construction contains virtually
all line bundles on this moduli space. The proof relies heavily
on the Grothendieck-Riemann-Roch theorem and the knowledge of the
singularities of the moduli space.

The theorem we will prove is the foolowing:

\pro{Main theorem}
Let $(X, H)$\ be any polarized algebraic sueface and $I\in\pic(X)$. Suppose
$H^2(X,\ZZ)$\ has no torsions, then there is an integer $N$\
depending on $(X,I,H)$\ such that for any $d\geq N$\ there is a homomorphism
$${\vphi}:\picxs\oplus \ZZ\lra \pic\bl\mhd\br\tz
$$
that has finite kernel and cokernel.
\endpro

\pro{Acknowledgment} The author would like to thank J. Koll\'ar
and K. Yoshioka for valuable discussion.
\endpro


\head Construction of line bundles on moduli spaces
\endhead

We begin with a general construction of line bundles on moduli space of
stable sheaves based on Knudsen and Mumford's determinant line bundle
construction. One way of getting line bundles on the moduli space $\m$\
is by forming (perfect) complexes on $\m\times X$\ out of the universal
sheaf $\f$\ and then take the determinant of its push-forward on $\m$.
This has been worked out by many people in various settings. It is
interesting to observe that we can also construct complexes on $\m$\ that
are push-forwards of complexes on $\m\times\xxx$\ constructed
based on $\pi_{12}\sta\f\otimes\pi_{13}\sta\f$, where $\pi_{12},\pi_{13}\mh
\m\times\xxx\to\m\times X$\ are projections. In this way, we recover the
previously known line bundles as well as a group of new line
bundles. We will work out the detail of this
construction in this section.

We begin with some word on the convention that will be used
throughout this paper. We let $X$\ be a fixed smooth algebraic surface
over complex numbers $\CC$\
and $H$\ a very ample line bundle on $X$. For any choice of $I\in\pic(X)$\
and $d\in\ZZ$($\cong H^4(X;\ZZ)$),
we form the (coarse) moduli scheme of $H$-stable
rank 2 sheaves $\e$\ satisfying $\det\e\cong I$\ and $c_2(\e)=d$,
which will be denoted by $\mhd$. (We will abbreviated it to $\md$\
when the choice of $H$\ is understood.)
We will denote by $\mbd$\ the (coarse) moduli scheme of $H$-semistable
sheaves satisfying the same restrain, modulo certain equivalence relation.
Following Gieseker [Gi1], $\md$\ is quasi-projective and $\mbd$\
is projective. Throughout this paper, we will
use Roman letters to denote vector bundles and use Calligraphic letters
to denote sheaves.
We denote by $A^iX$\ (resp. $A_iX$) the Chow cohomology (resp. homology)
group of $X$. For any sheaf, we denote by $\tilde c_i$\ the Chern class
taking value in the Chow cohomology group and by $c_i$\ the Chern class
taking value in the (ordinary) cohomology group. We denote by $\ch$\ the Chern
character taking value in $A\sta\lowq$, the Chow cohomology group with
rational coefficient.
When we assign a letter, say $\pi$, to denote projection from a
product space to its factor(s), we will use $\pi_i$\ to denote projection
to the $i$-th copy and use $\pi_{ij}$\ to denote projection to the
produce of $i$-th and $j$-th copies. Sometimes we will simply use, say $\pi_X$,
to denote projection to the factor $X$.

Now we construct line bundles on $\md$.
For the moment, we assume $\md$\ admits a universal family, namely a
sheaf $\f$\ on $\md\time X$\ flat over $\md$\ whose restriction to fiber
$\{w\}\time X$\ is isomorphic to the sheaf represented by $w\in\md$.
Note that the choice of universal family is not unique. Any two such families
differ by tensoring a pullback invertible sheaf on $\md$.
Let $K$\ be the Grothendieck's
group of $X\time X$\ and let $\tilk$\ be the kernel of $\xi\mh K\to\ZZ$,
$$\xi(C)=\xx(C\tor \pia\sta\e\tor\pib\sta\e)
\tag 1.1
$$
for some $\e\in\md$\ and $\pi_i\mh\xxx\to X$\ projection.
(Here and later, we will use $\tor$\ to
denote the tensor product of complexes.) We let $p_{\,\cdot}$\
be the projection from $\md\time\xxx\to\md$\ to its factor(s) with
self-explanatory index. Then we can define a homomorphism
$$\tilk\lra \picmd
\tag 1.2
$$
as follows: Since the universal family $\f$\ is flat over $\md$, $\f$\ is a
perfect complex on $\md\time X$. Then for any element in $\tilk$\
represented by a complex $C$\ on $\xxx$, we can form the complex
$\pbc\sta C\tor \pab\sta\f\tor\pbc\sta\f$\ and the push forward complex
$$\bl\pa\br_!\bl\pbc\sta C\tor \pab\sta\f\tor\pbc\sta\f\br.
$$
Following [KM], we can define the determinant line bundle
$$\det\Bl\bl\pa\br_!\bl\pbc\sta C\tor \pab\sta\f\tor\pbc\sta\f\br\Br.
\tag 1.3
$$
The homomorphism (1.2) is the one sending $C\in\tilk$\ to the line bundle
(1.3). Observe that if we replace $\f$\ by a different universal family, say
$\f\pri=\f\otimes\pia\sta\OO(L)$, $\pia\mh\md\time X\to\md$, then
$$\det\Bl\bl\pa\br_!\bl\pbc\sta C\tor \pab\sta\f\pri\tor\pbc\sta\f\pri\br\Br
=L^{\otimes(2\xi(C))} \otimes
\det\Bl\bl\pa\br_!\bl\pbc\sta C\tor \pab\sta\f\tor\pbc\sta\f\br\Br,
$$
where $\xi(C)$\ is as in (1.1). Hence
(1.2) is independent of the choice of the universal family since
$C\in \tilk$\ and then the homomorphism (1.3) is well-defined.

We will see shortly that the homomorphism (1.2) is far from injective.
In the following we will
use the Grothendieck-Riemann-Roch theorem (abbreviate
GRR theorem) to give an explicit description of the image of (1.2).
Since $\pa\mh\md\times\xxx\to\md$\ is smooth,
$$\ch\Bl\bl\pa\br_!\bl\picff\br\Br
=p_{1\ast}\Bl\ch\bl\picff\br\cdot\pbc\sta\tdxx\Br,
\tag 1.4
$$
by [Fu, Example 18.3.10]. By choosing $C\in\tilk$, the summand
in $A^0\bl\md\br\lowq$\ of (1.4) is trivial and the summand in
$A^1\bl\md\br\lowq$\ reads
$$\align
&\co\Bl\bl\pa\br_!\bl\picff\br\Br\\
=&\Bigl[ p_{1\ast}\bl\ch(\picff)\cdot \pbc\sta\tdxx\br\Bigr]\lowo.
\tag 1.5
\endalign
$$
We will use $\bigl[\cdot\bigr]_{[i]}$\ to denote its component in $A^i$.
A simple calculation shows that it is of the form
$$\align
\Bigl[ p&_{1\ast}\Bl \bl\pab\sta c_2(\f)+\pac\sta c_2(\f)\br\cdot\alp_1
+\bl\pab\sta c_3(\f)+\pac\sta c_3(\f)\br\cdot\alp_2\\
&+\pab\sta c_2(\f)\cdot\pac\sta c_2(\f)\cdot\alp_3
+\bl\pab\sta c_2(\f)\cdot \pac\sta c_3(\f)+
\pac\sta c_2(\f)\cdot \pab\sta c_3(\f)\br\cdot\alp_4\Br\Bigr]\lowo,
\endalign
$$
for some $\alp_{\,\cdot}\in A\sta(\xxx)\lowq$\ based on type consideration
and
$$\ch(\f)\equiv 2-\ct(\f)+{1\over 2}\tilde c_3(\f) \mod\co(\f).
$$
Hence (1.5), which is the Chern class of (1.3),
is contained in the linear span of the images
$$p_{M\ast}\bl \ct(\f)\cdot p_X\sta(\cdot)\br
: A^1(X)\lra A^1\bl\md\br\lowq,
\tag 1.6
$$
$$p_{M\ast}\bl \tilde c_3(\f)\cdot p_X\sta(\cdot)\br
: A^0(X)\lra A^1\bl\md\br\lowq
\tag 1.7
$$
($p_M, p_X$\ are projections of $\md\times X$) and
$$p_{1\ast}\bl\pab\sta\ctf\cdot\pac\sta\ctf\cdot\pbc\sta(\cdot)\br
: A^1(\xxx)\lra A^1\bl\md\br\lowq.
\tag 1.8
$$
It is easy to check that the image
of (1.6) is contained in the image of (1.8) when $4d>I^2$, which we assume
in the sequel. Thus the Chern class of (1.3) is contained in the image of
$$A^1\bl\xxx\br\lowq\oplus\QQ\lra
A^1\bl\md\br\lowq
\tag 1.9
$$
that is a direct sum of (1.8) and (1.7).

We remark that since $\md$\ is singular in general, the homomorphism
given by the first Chern class
$$\co:\picmd\lra A^1\bl\md\br\lowq
$$
may have infinite kernel. However, (1.9) suggests that very likely
there should be a homomorphism
$$\vphi : \picxs\oplus \ZZ\lra\picmd
\tag 1.10
$$
compatible to (1.9) that will capture all (up to finite index) line bundles
constructed in (1.3). Here $\sig\mh \xxx\to\xxx$\ is the map exchanging
factor and $\picxs$\ consists of line bundles invariant under $\sig$. We use
$\picxs$\ instead of $\pic\bl\xxx\br$\ in (1.10) since the map (1.2) so
constructed is $\sig$\ invariant.

We now construct this homomorphism explicitly. We first let
$$\psi:\picxs\oplus\ZZ\lra \tilk\tz
$$
be the map sending $L\in\picxs$\ to
$$\bl\OO_{\xxx}-\OO(L)\br+\half\bl\OO_{\xxx}-\OO(L)\br^2-{1\over 3}
\bl\OO_{\xxx}-\OO(L)\br^3+\alp \CC_{p_0}\in\tilk\otimes_{\ZZ}\QQ
$$
and sending $m\in\ZZ$\ to
$$C_0=m\OO_{X\time\{x_0\}}+\alp\CC_{p_0}\in\tilk\otimes_{\ZZ}\QQ,
$$
where $p_0\in\xxx$\ is a fixed point and $\alp\in \ZZ[{1\over 12}]$\ is
chosen so that $\psi(L)$\ and $\psi(m)$\ lie in $\tilk\tz$.
($\alp\in\ZZ[{1\over 12}]$\
follows from $L$\ being $\sigma$-invariant.)
Note that this choice gives
$$\ch\bl\psi(L)\br=\co(L)+\alp[p_o]\dual.
$$
We then define $\vphi(\alpha)$,
$\alpha\in\picxs\oplus\ZZ$, to be the $\QQ$-line bundle
$$\det\Bl\bl\pa\br_!\bl\pbc\sta\psi(L)\tor\ff\br\Br\in\picmd\tz.
$$
This way we have defined
$$\vphi:\picxs\oplus\ZZ\lra \picmd\tz.
\tag 1.11
$$
$\vphi$\ is a homomorphism, which can be proved similar to that
on page 422 of [Li1] despite $\psi$\ is not necessarily a homomorphism.
(One hint that this is true comes from the identity
$\ch\bl\psi(L)\br=\co(L)+\alp[p_o]\dual$.)
We shall omit the proof here since it is rather routine.

In the remainder of this section, we shall investigate how to define
this homomorphism when there are no universal family
on $\md\time X$\ and under what condition does it extend to
$$\vphib:\picxs\oplus\ZZ\lra \picmdb\tz.
$$
We will accomplish this by taking $\md$\ as a GIT quotient of a
Quot-scheme and then apply the descent theory.
Following [Gi1, Ma], there is a Quot-scheme $\qu$\ with a linear
$G=SL(N)$\ action such that if we denote
by $\qu^s$\ (resp. $\qu\ss$) the open subset of $G$-stable (resp.
$G$-semistable) points, then
$$\qu^s/G=\md\quad\and\quad
\qu\ss\quot G=\mdb.
$$
Let $\e$\ be a universal quotient family on $\qu\ss\times X$.
Then we can define determinant line bundles on $\qu\ss$\ as
we did for $\md$\ by using the family $\e$\ and
complexes in $\tilk$. Let $\tilde p_{\,\cdot}$\ be
projections from $\qu\ss\time\xxx$\ to its factor(s) with self-explanatory
index. Then for any $C\in\tilk$, we get a line bundle on $\qu\ss$:
$$L(C)=\det\Bl\bl\tpa\br_!\bl\tpbc\sta C\tor\tpab\sta\e\tor\tpac\sta\e\br\Br.
\tag 1.12
$$
$L(C)$\ admits a canonical $G$-linearization since the $G$-action on
$\qu\ss$\ lifts to a $G$-action on $\e$.
By descent theory, $L(C)$\ descends to $\md$\ if for any $w\in\qu^s$\
the stabilizer $\stab_w\sub G$\
of $w$\ acts trivially on the fiber of $L(C)$\ over $w$. For the same
reason, it descends to
$\mdb$\ if similar condition holds for all $w\in\qu\ss$\ having closed
orbits $G\cdot w$.

\pro{Definition 1.1}
An ample line bundle $H$\ is said to be $(I,d)$-generic
if whenever $J\in\pic(X)$\
satisfies $\bl2c_1(J)-c_1(I)\br\cdot c_1(H)=0$, then
$-4c_1(J)^2>4d- c_1(I)^2$\ unless $ c_1(J)$\ is torsion.
\endpro

Note that if $H$\ is $(I,d)$-generic, then any strictly semistable sheaf
$\f\in\mdb$\ is S-equivalent to $\f_1\oplus \f_2$\ with $2c_1(\f_1)\equiv
c_1(I)$\ modulo torsions.

\pro{Proposition 1.2}
For any $C\in\tilk$, the line bundle $L(C)$\ descends to $\md$. Further
if $H$\ is $(I,d)$-generic, then it descends to $\mdb$.
\endpro

\proof
For stable point $w\in\qu^s$, the stabilizer $\stab_w$\ is $\ZZ/N\ZZ$\
acting on $\e\otime k(w)$\ via multiplication by an N-th root of unity,
say $\epsilon$. Here $\e\otime k(w)$\ is the restriction of $\e$\ to the
fiber $\{w\}\time X$. Then its induced action on $L(C)_w$\ (which is
the restriction to the fiber over $w$) is a multiplication
by $\epsilon^m$, where $m$\ is given by (1.1) that is 0
since $C\in \tilk$. Therefore by descent theory $L(C)$\ descends to $\md$.

We now study descent problem on $\mdb$. Let
$w\in\qu\ss$\ be a strictly semistable point with closed orbit $G\cdot w$.
Then $\e\otimes k(w)$\ is a direct sum of
rank one sheaves, say $\f_1\oplus \f_2$. If we further assume $H$\ is
$(I,d)$-generic, then $c_1(\f_1)$\ and $c_1(\f_2)$\ differ by a
torsion element. In this case, it is
straight forward to check, as did on page 426 of [Li1], that $\stab_w$\
acts trivially on $L(C)_w$. Therefore, $L(C)$\ descends to
$\mdb$\ as claimed.
\qed

\pro{Remark} When $H$\ is not $(I,d)$-generic, one can determine
explicitly which line bundles do descend by using the descent argument,
if one knows all strictly semistable sheaves.
\endpro

Knowing that (1.12) always descends to $\md$, we can define a
homomorphism
$$\vphi:\picxs\oplus\ZZ\lra\picmd\tz
$$
that sends $\alp$\ in $\picxs\oplus\ZZ$\ to the descent of the
determinant line bundle associated to the complex
$\psi(\alp)$. Similarly, if $H$\ is $(I,d)$-generic then the
same construction yields a homomorphism
$$\vphib:\picxs\oplus\ZZ\lra \picmdb\tz.
\tag 1.13
$$

\head 2. Injectivity of the homomorphism $\vphi$
\endhead

In this section, we will show that for sufficiently large $d$, the homomorphism
$$\vphi:\picxs\oplus \ZZ\lra\picmd\tz
$$
has finite kernel. The method we will use is to construct a
variety $W$\ and a morphism
$$g\mh W\to \md
\tag 2.1
$$
and show that the composition
$$g\sta\circ\vphi:\picxs\oplus\ZZ\lra\pic(W)\tz
$$
has finite kernel, using GRR theorem.

The variety $W$\ we shall use is the one parameterizing non-locally free
sheaves that are kernels of $\e\to \CC_x\oplus\CC_y$, where $\e$\ is
fixed and $x,y\in X$\ are arbitrary closed points. $W$\ will be birational
to a $\bold P^1\times \bold P^1$\ bundle over $\xxx$.
Let $\e$\ be a rank two locally free sheaf on $X$\ and let $\pe$\ be
the associated projective bundle with projection
$$\pi : \pe\to X
$$
(we adopt the convention in [Ha, p162] that $\pi_{\ast}\bl \oo\br=\e$)
and let $W$\ be the blowing-up of $\pe\times\pe$\ along the diagonal. We first
construct a sheaf on $W\times X$\ flat over $W$. Let $q_i$\ be the
projection
$$q_i: W\lra\pe\times\pe\mapright{\text{pr}_i}\pe,\quad i=1,2
$$
and let
$$\iota_i:\text{id}\times\bl\pi\circ q_i\br: W\lra W\times X,\quad i=1,2.
\tag 2.2
$$
$\iota_i$\ is a closed immersion. We let
$$\beta:\pi_X\sta \e\lra \iota_{1\ast}q_1\sta \oo\oplus \iota_{2\ast}q_2\sta\oo
$$
be the homomorphism induced by
$$\pi_X\sta \e\lra \pix\sta \e\otimes\OO_{\iota_i(W)}=
\iota_{i\ast}q_i\sta\pi\sta \e \lra\iota_{i\ast}q_i\sta\oo,
$$
where $\pix\mh W\time X\to X$\ is the second projection.
$\beta$\ is surjective away from $S\time X$, where $S$\ is the exceptional
divisor of $W$. We let $\g$\ be the kernel of $\beta$\ and let
$\Omega$\ be the cokernel of
$$\g\lra \pi_X\sta\e.
$$
In this way, we obtain two exact sequences:
$$0\lra\g\lra\pix\sta \e\lra \Omega\lra 0
\tag 2.3
$$
and
$$0\lra\Omega\lra\iota_{1\ast}q_1\sta \oo\oplus \iota_{2\ast}q_2\sta\oo
\lra\iota_{1\ast}q_1\sta \oo_{|S\times X}\lra 0
\tag 2.4
$$
since $\coker(\beta)$\ is isomorphic to the last term in (2.4).

We claim that $\g$\ is a family of torsion free sheaves flat over $W$. Because
of the exact sequence (2.3), it suffices to show that $\Omega$\ is flat over
$W$\ and
$$\operatorname{Tor}\bl\Omega,\OO_{\{w\}\times X}\br=0
$$
for all closed $w\in W$. But this is clear from the exact sequence
(2.4) because $\iota_i(W)\sub W\times X$\ is smooth over $W$,
$q_i\sta\oo$\ is an invertible sheaf on
$W$\ and $\iota_i(W)\cap (S\time X)$\ is a divisor in $\iota_i(W)$.
This proves the claim.

Following the discussion in \S 1, we can form a homomorphism
$$\vphi_W:\picxs\oplus\ZZ\lra\pic(W)\tz
\tag 2.5
$$
using the determinant line bundle construction outlined before (1.11)
using the sheaf $\g$\ over $W\times X$\ and complexes
$$\psi:\picxs\oplus\ZZ\lra \tilk\tz.
$$
We will investigate the injectivity of this homomorphism by looking at
their first Chern classes, using GRR. From the formula (1.5) and the
one after, for $L\in\picxs$
$$\co\bl\vphi_W(L)\br=\Bigl[ p_{1\ast}\Bl\ch\bl\plgg\br\cdot\pbc\sta
\tdxx\Br\Bigr]_{[1]},
$$
where $p_{\,\cdot}$\ are projections from $W\time \xxx$\ to its factor(s)
with self-explanatory index. Since fibers of $p_1$\ are $\xxx$, it suffices
to determine
$$\align
\Bigl[\ch\bl&\plgg\br\cdot\pbc\sta\tdxx\Bigr]_{[5]}\\
&=\Bigl[\pbc\sta\bl\co(L)+\alp[p_0]\dual\br\cdot
\pab\sta\ch(\g)\cdot\pac\sta\ch(\g)\cdot
\pb\sta\tdx\cdot\pc\sta\tdx\Bigr]_{[5]}
\tag 2.6
\endalign
$$
modulo the kernel of $p_{1\ast}$. We first determine $\ch(\g)$. By using
(2.3) and (2.4) we get
$$\ch(\g)=\ch\bl\pix\sta \e\br-\ch\bl\iota_{1\ast}q_1\sta\oo\br
-\ch\bl\iota_{2\ast}q_2\sta\oo\br+
\ch\bl\iota_{1\ast}q_1\sta\oo_{|S\times X}\br.
\tag 2.7
$$
($\pi_{\,\cdot}$\ is a projection from $W\times X$\ to its factor.)
Applying GRR theorem to the proper morphism $\iota_1\mh W\to W\time
X$, we obtain
$$\align
\ch\bl\iota_{1\ast}q_1\sta\oo\br&=
\iota_{1\ast}\bl\ch(q_1\sta\oo)\br\cdot\td(W)\br\cdot\td(W\time X)\umo\\
&=\iota_{1\ast}q_1\sta e^{\co(\oo)}\cdot \pi_W\sta\td(W)
\cdot\td(W\time X)\umo\\
&=\iota_{1\ast}q_1\sta e^{\co(\oo)}\cdot \pix\sta\tdx\umo.
\endalign
$$
Here the second identity holds because
$\iota_1(W)$\ is a graph of $\pi\circ q_1\mh W\to X$. Similarly, we
have
$$\ch\bl\iota_{2\ast}q_2\sta\oo\br
=\iota_{2\ast}q_2\sta e^{\co(\oo)}\cdot \pix\sta\tdx\umo.
$$

In the following, we will choose a split $\e$\ (i.e. $\e$\ is a
direct sum of invertible sheaves) to simplify our calculation.
We assume $\det \e=I$. In this case, $\pe$\ is isomorphic to $\bold P^1
\time X$, $\co(\oo)$\ is dual to a section of $\pe\to X$\ and
$\co\bl\oo\br^2=0$.

We now reorganize terms in the expansion
of (2.6) in terms of (2.7) for split $\e$. One term in the expansion of
(2.6) is
$$\align
\Bigl[\pbc\sta\bl\co(L)+\alp[p_0]\dual
\br\cdot(\pab\circ\pix)\sta & \ch(\e)\cdot
(\pac\circ\pix)\sta\ch(\e)\Bigr]\lowf\\
&=\Bigl[\pbc\sta\bl\co(L)+\alp[p_0]\dual
\br\cdot\pb\sta\ch(\e)\cdot\pc\sta\ch(\e)\Bigr]\lowf
\endalign
$$
which is contained in the kernel of $p_{1\ast}$, following the
vanishing of (1.1).
(Recall $p_{\,\cdot}$\ is designated to projections from $W\times\xxx$\
to its factor, $\pi_{\,\cdot}$\ is projection from $W\times X$\ to its
factor and $q_{\,\cdot}$\ is projection from $W$\ to $\pe$.)
The next terms are
$$\align
\Bigl[\pbc\sta\co(L)\cdot\pab\sta\bl&\iota_{i\ast}q_i\sta e^{\co(\oo)}
\cdot\pix\sta\tdx\umo\br\\
&\cdot\pac\sta\bl\iota_{j\ast}q_j\sta e^{\co(\oo)}\cdot\pix\sta\tdx\umo\br
\cdot\pb\sta\tdx\cdot\pc\sta\tdx\Bigr]\lowf,
\endalign
$$
where $(i,j)$\ runs through pairs $(1,1)$, $(1,2)$, $(2,1)$\ and $(2,2)$.
After simplification, their images under $p_{1\ast}$\ are
$$p_{1\ast}\Bl\pbc\sta\co(L)\cdot\pab\sta\bigl[\iota_i(W)\bigr]\dual\cdot
\pac\sta\bigl[\iota_j(W)\bigr]\dual\Br
=\cases \qxx\sta\co(L),&i\ne j;\\
\qxx\sta \pr_1\sta\bl\co(L_{|\Delta})\br,& i,j=1;\\
\qxx\sta \pr_2\sta\bl\co(L_{|\Delta})\br,& i,j=2,
\endcases
$$
since $L$\ is invariant under permutation $\sigma$. Here
$\qxx\mh W\to \xxx$\ is the projection, $\Delta\sub\xxx$\ is
the diagonal, $\co(L_{|\Delta})$\ is in $A^1X$\ by the isomorphism
$\Delta=X$\ and $\pi_{\,\cdot}\mh \xxx\to X$\ is projection.
Note that $\qxx\sta\co(L)$\ appears in the expression.

We need to take care of the remainder terms. We first look at
another series of terms
$$\align
-\Bigl[\pbc\sta\co(L)\cdot (p_{1i}\circ\pix)\sta & \ch(\e)\cdot p_{1\bar i}
\sta\bl \iota_{1\ast} q_1\sta e^{\co(\oo)}\cdot\pix\sta\tdx\umo\\
&+\iota_{2\ast} q_2\sta e^{\co(\oo)}\cdot\pix\sta\tdx\umo\br\cdot\pb\sta\tdx
\cdot\pc\sta\tdx\Bigr]\lowf,
\tag 2.8
\endalign
$$
where $(i,\bar i)$\ is either $(2,3)$\ or $(3,2)$. For simplicity, we
let $\alp_2$\ and $\alp_3\in A\sta(\xxx)$\ be such that
$$-\pbc\sta\alp_i=\pbc\sta\co(L)\cdot
p_i\sta\ch(\e)\cdot p_i\sta\tdx.
$$
For our purpose, it is sufficient to look at their
components in $\qxx\sta A^1(\xxx)$\ under the decomposition
$$A^1(W)=\qxx\sta A^1(\xxx)\oplus \ZZ[R_1]\dual\oplus \ZZ[R_2]\dual\oplus
\ZZ[S]\dual,
$$
where $R_1,R_2\sub W$\ are birational to
$\PP^1\times\{\text{pt}\}\times\xxx$\ and
$\{\text{pt}\}\times\PP^1\times\xxx$\ respectively, $S\sub W$\ is the
exceptional divisor of $W\to\pe\times\pe$. (Recall $W$\ is a blowing up
of $\PP^1\times\PP^1\times\xxx$\ since $\e$\ is split.)
Since $\iota_i(W)$\ is
a graph of $W\to X$, the component in $\qxx\sta A^1(\xxx)$\ of the
image under $p_{1\ast}$\ of (2.8) are
$$\Bigl[ p_{1\ast}\Bl\pbc\sta\alp_2\cdot\pac\sta
\bl[\iota_1(W)]\dual+[\iota_2(W)]\dual\br\Br\Bigr]_{[1]}
\tag 2.9
$$
and
$$\Bigl[ p_{1\ast}\Bl\pbc\sta\alp_3\cdot\pab\sta
\bl[\iota_1(W)]\dual+[\iota_2(W)]\dual\br\Br\Bigr]_{[1]}.
\tag 2.10
$$
To simplify these, we proceed as follows: Consider the projection $f\mh
W\to\xxx$\ (from $W\to\PP(\e)\times\PP(\e)$) and the projection
$$P=(f,\text{id}_X,\text{id}_X): W\times\xxx\lra\xxx\times\xxx.
$$
Then by definition of $\iota_1$\ and $\iota_2$, we have
$$\pab\sta\bl[\iota_j(W)]\dual\br=
P\sta\circ\pi_{j3}\sta\bl[\Delta]\dual\br
$$
and
$$\pac\sta\bl[\iota_j(W)]\dual\br=
P\sta\circ\pi_{j4}\sta\bl[\Delta]\dual\br,
$$
where $\pi_{ij}\mh\xxx\times\xxx\to\xxx$\ is a projection. Then the
component in $A^1\bl\xxx\br$\ of (2.9) is the pull back (under
$q_{XX}\sta$) of
$$\Bigl[\pi_{12\ast}\circ P_{\ast}
\Bl[\Sigma\times\xxx]\dual\cdot
\bl P\circ\pi_{34}\br\sta(\alp_2)\cdot
P\sta\bl\pi_{14}\sta\bl[\Delta]\dual\br+\pi_{24}\sta\bl[\Delta]\dual)\br\Br
\Bigr]_{[1]},
$$
where $\Sigma$\ is a birational section of $W\to\xxx$. By projection
formula, it is
$$\align
&\Bigl[\pi_{12\ast}\Bl\bl P_{\ast}[\Sigma\times\xxx]\dual\br\cdot
\pi_{34}\sta(\alp_2)\cdot
\bl\pi_{14}\sta([\Delta]\dual)+\pi_{24}\sta([\Delta]\dual)\br\Br
\Bigr]_{[1]}\\
=&\Bigl[\pi_{12\ast}\Bl\pi_{34}\sta(\alp_2)\cdot
\bl\pi_{14}\sta([\Delta]\dual)+\pi_{24}\sta([\Delta]\dual)\br\Br
\Bigr]_{[1]}\\
=&\Bigl[\pr_1\sta\pr_{1\ast}(\alp_2)+\pr_2\sta\pr_{1\ast}(\alp_2)
\Bigr]_{[1]}.
\endalign
$$
Similarly, (2.10) is
$$\Bigl[\pr_1\sta\pr_{2\ast}(\alp_3)+\pr_2\sta\pr_{2\ast}(\alp_3)\Br
\Bigr]_{[1]}.
$$
Note that $\sigma\sta(\alp_2)=\sigma_3$, where $\sigma$\ is the
permutation. Hence both are
$$\pr_1\sta(\alpha)+\pr_2\sta(\alpha),\quad \alp=\pr_{1\ast}(\alp_2).
$$
The image under $p_{1\ast}$\ of the
last set of terms in (2.6) contains factor
$\ch\bl\iota_{i\ast} q_i\sta\oo_{|S\times X}\br$\
whose components in $\qxx\sta A^1(\xxx)$\ are trivial.
Put them together, we see that the component in $\qxx\sta A^1(\xxx)$\
of $\co\bl\vphi_W(L)\br$\ is the pull-back of
$$2\co(L)+\pr_1\sta\bl\co(L_{|\Delta})\br+
\pr_2\sta\bl\co(L_{|\Delta})\br+2\pr_1\sta(\alpha)
+2\pr_2\sta(\alpha)\in A^1(\xxx).
\tag 2.11
$$

Now we look at the component in $\qxx\sta A^1(\xxx)$\ of
$\co\bl\vphi_W(1)\br$. By replacing $\co(L)$\ in (2.6) with
$\bigl[ X\times\{x_0\}+\{x_0\}\times X\bigr]\dual$\ and going through
the same line of analysis, we see that there is a
$\beta\in A^1(X)$\ depending only on $X$\ and
$I$\ such that the component in $\qxx\sta A^1(\xxx)$\ of
$\co\bl\vphi_W(1)\br$\ is
$$\qxx\sta\bl\pr_1\sta\beta+\pr_2\sta\beta\br.
\tag 2.12
$$
Therefore, if $(L,n)\in\picxs\oplus\ZZ$\ makes
$$\co\bl\vphi\bl(L,n)\br\br=0,
$$
then $\co(L)$\ must be of the form $\pr_1\sta c+\pr\sta_2 c$\ for some
$c \in A^1(X)$\ satisfying the identity
$$ F(c)+c_2(\e)c+n\beta=0,
\tag 2.13
$$
where $F\mh A^1X\to A^1X$\ is a map depending on $X$\ and $I$\ but not
$\e$, from (2.11) and (2.12). (The exact form of (2.13) can be
derived easily from (2.11) and (2.12) but we
choose not to in part to emphasize the term $c_2(\e)\,c$\ that is
crucial in our later discussion.)

Now we state and prove the main proposition of this section:

\pro{Proposition 2.1}
For any choice of $(I,H)$, there is an $N$\ such that whenever
$d\geq N$, then the homomorphism
$$\vphi:\picxs\oplus\ZZ\lra\picmd\tz
$$
has finite kernel.
\endpro

\proof
We let $\e_1$\ be $\OO\oplus\OO(I)$\ and let $\e_2$\ be
$\OO(H)\oplus\OO(I\otimes H\umo)$. (When $H^2=H\cdot I$\ then we
let $\e_2$\ be $\OO(H^{-1})\oplus\OO(I\otimes H)$.)
Let $W_i$\ be the variety constructed at the beginning of this section
that is the blowing-up of $\PP(\e_i)\times\PP(\e_i)$\ and
let $\g_i$\ be the sheaf on $W_i\times X$\ constructed in (2.3).
We claim that if
$$(L,n)\in\picxs\oplus\ZZ
$$
lies in the kernel of $\vphi_{W_i}$\ in (2.5) for $i=1,2$, then
$n=0$\ and $L$\ is torsion. Indeed, by the previous argument,
$\vphi_{W_i}\bl(L,n)\br=0$\ implies that $\co(L)=\pr_1\sta c+\pr_2\sta c$\
with $c\in A^1X$\ satisfying two identities that are (2.8) with $c_2(\e)$\
replaced by $c_2(\e_1)=0$\ and by $c_2(\e_2)\ne0$. Clearly, this is
possible only if $c$\ is a torsion.
Once we know that $c$\ is torsion, then $\vphi_W\bl(L,n)\br=0$\
force $n=0$, which can be checked directly. This establishes the claim.

Now we prove the proposition. Let $(L,n)$\ be any element in the kernel
of $\vphi$. It suffices to show that $\vphi_{W_i}\bl(L,n)\br=0$\ for both $i$.
Let $\e$\ be either $\e_1$\ or $\e_2$. Then there
is a constant $N\geq 0$\ such that for any $d\geq N$\ there is a
subsheaf $\f\sub\e$\ such that $\det\f=I$\ and $c_2(\f)=d$\ satisfying the
following property: $\f$\ admits a deformation, say $\f_t$, whose
general members are locally free and $\mu$-stable.
This can be proved as follows: First for large $d$, we can choose
$\f\sub\e$\ so that the traceless part $\text{Ext}^2(\f,\f)^0=0$. Then
we can deform $\f$\ to locally free sheaves. When $d$\ is large enough
and the support of $\e/\f$\ is in general position, then the argument on
page 158 of [Gi2] shows that any locally free deformation of $\f$\ are
$\mu$-stable. This proves the existence of such deformation. Now let
$0\in T$\ be a smooth curve and $\f_T$\ a deformation of $\f_0=\f$\ just
mentioned. Without loss of generality, we assume $\f_T$\ is locally free
over $T-0$. Let $\Lambda\sub X$\ be the finite set over which $\f_0$\ is
not locally free.

We now construct the corresponding family of varieties
$W_T$\ associated to $\f_T$. Let $\PP(\f_T)$\ be the projective bundle
of $\f_T$\ over $T\times X-\{0\}\times \Lambda$\ and let $W_T$\
be the blowing-up of the diagonal of
$$\PP(\f_T)\times_T\PP(\f_T).
$$
Similar to (2.3), we can define a sheaf $\g_T$\ on $W_T\times X$\ that is
flat over $W_T$. Note that since $\f_t$\ is $\mu$-stable for $t\ne0$,
the family $\g_{T\sta}$, $T\sta=T-0$,
that is the restriction of $\g_T$\ to $W_T-W_0$\
($W_0$\ is the fiber over $0\in T$) is a family of $H$-stable sheaves
with determinant $I$\ and second Chern class $d$. Hence, it induces a
morphism
$$\mu : W_T-W_0\lra \md.
$$
Also, by applying the determinant line bundle construction to the family
$\g_T$, we can form a homomorphism
$$\vphi_{W_T}: \picxs\oplus\ZZ\lra\pic\bl W_T\br\tz
$$
similar to (2.5).
Since this construction commutes with base change,
$$\vphi_{W_T}\bl(L,n)\br_{|W_T-W_0}\cong \mu\sta\vphi\bl(L,n)\br.
$$
In particular, since $\vphi\bl(L,n)\br$\ is trivial,
$\vphi_{W_T}\bl(L,n)\br_{|W_T-W_0}$\ is trivial and consequently
$\vphi_{W_T}\bl(L,n)\br_{|W_0}$\ is trivial.

Assume $\f^{\vee\vee}=\e$\ is $\e_1$\ we begin with. Then
$W_0$\ is canonically isomorphic to $W_1-A$\ for some closed $A\sub W_1$\
of codimension two and further,
$$\vphi_{W_T}\bl(L,n)\br_{|W_0}\cong \vphi_{W_1}\bl(L,n)\br_{|W_1-A}.
$$
Hence for $(L,n)\in\text{ker}\,\bl\vphi\br$, $\vphi_{W_1}\bl(L,n)\br$\
is trivial as well, since $\codim(A)=2$\ and $W_1$\ is smooth.
Similarly, if we choose $\f$\ so that $\f^{\vee\vee}=\e_2$, then we must
have $\vphi_{W_2}\bl(L,n)\br=0$.
This shows that $\vphi_{W_i}\bl(L,n)\br$\ are trivial for $i=1,2$\ and
henceforth establishing the proposition.
\qed

\head 3. Surjectivity of $\vphi$
\endhead

In this section, we shall show that when $H^2(X,\ZZ)$\ has no torsion, then
for sufficiently large $d$\ the so constructed homomorphism
$$\vphi:\picxs\oplus\ZZ\lra\picmd\tz
$$
has finite cokernel. The tactic we will use is to first show that
under certain circumstances the restriction homomorphism
$$\pic\bl\mdb\br\lra\picmd
$$
is surjective, by studying extension
problem. We will then use the knowledge of the first two Betti numbers
of $\md$\ to determine the Picard group of $\mdb$, up to finite index.

We first look at the problem of extending line bundle on $\md$\ to
$\mdb$. This in general is rather tricky due to the singularities of $\mdb$.
In out case, we can give an affirmative answer to this problem because of
our knowledge of singularities of $\mhdb$. To this end, we will
quickly review some relevant facts about the singularities of the
moduli space $\mdb$.
Since the discussion of $\picmdb$\ will simplify when the polarization
$H$\ is $(I,d)$-generic, we will work with a set of polarizations
simultaneously. (Note that except when $\pic(X)/\pic^0(X)$\ has
rank one, no polarization is $(I,d)$\ generic for all $d$.)
In the following, we will use $\mhd$\ and
$\mhdb$\ to denote the moduli schemes of $H$-stable and $H$-semistable
sheaves respectively. Also, we will speak freely of moduli scheme
for $\QQ$-ample line bundle $H$\ since $\mhd$\ and $\mhdb$\ only
depend on the ray $\QQ^+\cdot c_1(H)\sub H^{1,1}(X,\QQ)$.
We let $\nqp$\ be the $\QQ$-cone in $H^{1,1}(X,\RR)\cap H^2(X,\QQ)$\
spanned by Chern classes of $\QQ$-ample line bundles. We say
a neighborhood $\cc\sub\nqp$\ is precompact if the closure of $\cc$\
in $H^{1,1}(X,\RR)\cap H^2(X,\QQ)$\ is still contained in $\nqp$.
Note that for any $c\in\nqp$, a Euclidean ball of sufficiently small
radius centered at $c$\ in $\nqp$, after fixing an Euclidean metric
on $H^{1,1}(X,\RR)$, is a precompact neighborhood of $c$. We will
use $\cc\subb\nqp$\ to mean $\cc$\ is precompact. Also, by abuse of
notation we will say a $\QQ$-line bundle $H\in \cc$\ if $c_1(H)\in\cc$.

We collect the relevant properties of $\mdb$\ in the following
lemmas:

\pro{Lemma 3.1}
For any precompact $\cc\subb\nqp$, there is an integer $N$\ depending
on $(X,I,\cc)$\ such that for any $d\geq N$\ and $\QQ$-ample
line bundle $H\in\cc$\ (i.e. $c_1(H)\in\cc$),
\roster
\item $\mhdb$\ is normal, irreducible and of dimension predicted by
Riemann-Roch theorem;
\item $\mhd$\ is a local complete intersection scheme and the
codimension of the singular locus of $\mhd$\ is at least 5;
\item Let $\tilh\in\cc$\ be another $\QQ$-ample line bundle with
$c_1(\tilh)\in\cc$, then there are closed subsets $U_H\sub\mhd$\ and
$U_{\tilh}\sub\mthd$\ having codimension at least 5 respectively
such that $\mhd-U_H$\ is canonically isomorphic to $\mthd-
U_{\tilh}$\ as schemes, by identifying two points that represent
isomorphic sheaves.
\endroster
\endpro

\proof
(1) and (2) are proved on page 10 of [Li3] and (3) is proved on page
458 of [Li1] by using [Qi2].
\qed

As was mentioned in \S 2, the moduli $\mhdb$\ is constructed
as GIT quotient of a Grothendieck's Quot-scheme, say $\qu_H(I,d)$,
by a reductive group $G$. (We take $G=PGL(m)$\ this time.)
The property we need from this construction is stated in the following
lemma.

\pro{Lemma 3.2}
For any precompact $\cc\subb\nqp$, there is an integer $N$\ depending
on $(X,I,\cc)$\ such that for any $d\geq N$\ and
$H\in\cc$, we can construct a Grothendieck's Quot-scheme
$\qu_H(I,d)$\ and a reductive group $G$\ acting linearly on
$\qu_H(I,d)$\ of which the following holds:
\roster
\item
Let $\qu_H(I,d)^s$\ and $\qu_H(I,d)\ss$\ be the open subsets of
$G$-stable and $G$-semistable (with respect to the given
linearization) points in $\qu_H(I,d)$. Then
$$\qu_H(I,d)^s/G=\mhd,\quad\qu_H(I,d)\ss\quot G=\mhdb
$$
and the projection sends $w\in\qu_H(I,d)\ss$\ to $\f_w\in
\mhdb$, where $\f_w$\ is the corresponding quotient sheaf of $w$;
\item
$\qu_H(I,d)\ss$\ is a local complete intersection scheme whose
singular locus has codimension at least 5 in $\qu_H(I,d)\ss$;
\item Let $w\in \qu_H(I,d)\ss$\ be any point with closed orbit $G
\cdot w\sub \qu_H(I,d)\ss$, then $\stab_w=\CC\sta$\ if $\f_w\cong
J_1\oplus J_2$\ with $J_1\ne J_2$, $\stab_w=PGL(2,\CC)$\ if
$\f_w\cong J\oplus J$\ and $\stab_w=\{1\}$\ otherwise;
\item Assume $H$\ is $(I,d)$-generic, then
$$\codim\bl\qu_H(I,d)\ss-\qu_H(I,d)^s,\qu_H(I,d)\ss\br\geq 8.
$$
\endroster
\endpro

\proof
(1) and (3) is proved by [Gi1] (see also [Ma]) and (2) is proved on
page 5 of [Li2] and (4) follows from [Qi2].
\qed

The last is the information about the first two Betti numbers of
the moduli space $\mhd$\ proved recently in [Li3].

\pro{Lemma 3.3}
For any precompact $\cc\subb\nqp$, there is an integer $N$\
depending on $(I,d,\cc)$\ such that for any $d\geq N$\ and $H\in\cc$,
$$H^1\bl\mhd;\QQ\br\cong H^1\bl X;\QQ\br
$$
and
$$H^2\bl\mhd;\QQ\br\cong H^2\bl X;\QQ\br\oplus
\wedge^2 H^1\bl X;\QQ\br\oplus H^0(X;\QQ).
$$
\endpro

We also quote a lemma that concerns extending line bundles
over local complete intersection singularities.

\pro{Lemma 3.4}
Let $W$\ be any quasi-projective scheme having only local
complete intersection singularity and $\codim\bl\text{Sing}\,W,W\br
\geq 4$. Then for any closed subset $\Lambda\sub W$\
of codimension at least two, the canonical map induced by
inclusion
$$\pic(W)\lra\pic(W-\Lambda)
$$
is an isomorphism.
\endpro

\proof See [SGA 2] on page 132.
\qed

Now we draw some easy consequence from these lemmas
concerning the Picard group of $\mhd$. In the following, we fix
an $I$, a $\cc\subb\nqp$\ and the integer $N$\ specified in the
previous lemmas.

\pro{Lemma 3.5}
For any $H_1,H_2\in\cc$\ and $d\geq N$, $\pic\bl \m_{H_1}(I,d)\br$\
is canonically isomorphic to $\pic\bl \m_{H_2}(I,d)\br$, compatible
to the birational map in Lemma 3.1.
\endpro

\proof
Let $U_{H_1}\sub \m_{H_1}(I,d)$\ and $U_{H_2}\sub\m_{H_2}
(I,d)$\ be the closed subset given in Lemma 3.1.
It suffices to show that the pull back homomorphism
$\pic\bl \m_{H_i}(I,d)\br\to \pic\bl \m_{H_i}(I,d)-U_{H_i}\br$\
is an isomorphism, since $\m_{H_1}(I,d)-U_{H_1}$\ is canonically
isomorphic to $\m_{H_2}(I,d)-U_{H_2}$.
But this follows from Lemma 3.4 and (2) of
Lemma 3.1. This proves the lemma.
\qed

\pro{Lemma 3.6}
Suppose $H^2(X,\ZZ)$\ is torsion free. Then for $d\geq N$\ and
$(I,d)$-generic $H\in\cc$, the pull back homomorphism
$$\mu: \pic\bl \mhdb \br\lra \pic\bl\mhd\br
$$
is an isomorphism.
\endpro

\proof
It is easy to see that $\mu$\ is injective. Indeed, assume $L$\
is a line bundle on $\mhdb$\ admitting a non-vanishing
section over $\mhd$. Since $\mhdb$\ is normal and
$$\codim\bl\mhdb-\mhd,\mhdb\br\geq 2,
$$
this section extends to a non-vanishing section over $\mhdb$.
Hence $L$\ is trivial itself.

Now we show that $\mu$\ is
surjective. Let $L$\ be a line bundle on $\mhd$. Then its pull back
$\pi\sta L$\ on $\qu_H(I,d)^s$\ admits a canonical
$G$-linearization, where $\pi\mh\qu_H(I,d)^s\to\mhd$. By
Lemma 3.4 and  Lemma 3.2, $\pi\sta L$\ extends to a
line bundle $\widetilde{L}$\ on $\qu_H(I,d)\ss$\ and the
$G$-linearization of $\pi\sta L$\ extends as well. The key
observation is that $\widetilde{L}$\ always descends to a line bundle on
$\mhdb$, under our assumption on $H^2(X;\ZZ)$.
To see this, we need to check that for any $w\in\qu_H(I,d)\ss$\
with closed orbit $G\cdot w$\ and non-trivial stabilizer
$\stab_w$, $\stab_w$\ acts trivially on the fiber of
$\widetilde{L}$\ over $w$. Indeed, since $\widetilde{L}$\ is a
line bundle, this action is given by a character
$$\xx_w:\stab_w\lra\CC\sta.
$$
By (3) of Lemma 3.2, $\stab_w$\ can only take two forms,
$\CC\sta$\ and $PGL(2)$, unless it is trivial. Since $PGL(2)$\
has no non-trivial character, we are left to show that $\xx_w$\
is trivial even when $\stab_w=\CC\sta$. Following Drezet and Narasimhan
[DN], we know that if $w_1,w_2$\ are two points contained in the same
connected component of
$$\bigl\{w\in\qu_H(I,d)\ss\mid G\cdot w\ \text{closed and}\
\stab_w\ne\{1\}\bigr\},
\tag 3.1
$$
then $\xx_{w_1}=0$\ if and only if $\xx_{w_2}=0$.
By (3) of Lemma 3.2, $G\cdot w$\ closed and $\stab_w
\ne\{1\}$\ implies $\f_w=\f_1\oplus \f_2$. However, since
$H$\ is $(I,d)$-generic, $c_1(\f_1)-c_1(\f_2)$\ is torsion, which is
zero since $H^2(X,\ZZ)$\ has no torsions.
Then we can deform $\f_1\oplus\f_2$\ within the subset of split
quotient sheaves in $\qu_H(I,d)\ss$, then within (3.1), to $\f\oplus\f$.
Therefore $w$\ and the point $\tilde w$\ associated to the quotient
sheaf $\f\oplus \f$\ lie in the same component of (3.1).
Therefore $\xx_w=0$\ because
$\stab_{\tilde w}=PGL(2)$\ has no non-trivial
characters. This proves that $\widetilde{L}$\ descends to
a line bundle on $\mhdb$, which is an extension of $L$.
This proves that $\pic\bl\mhdb\br\to\pic\bl\mhd\br$\ is surjective.
\qed

As a consequence of this lemma, we see that for large $d$, the
divisor $\Lambda\sub\mhdb$\ consisting of non-locally free sheaves
is Cartier when $H$\ is $(I,d)$-generic and $H^2(X,\ZZ)$\ is torsion free.
In the following, we shall show that the condition $H^2(X;\ZZ)$\ has no
torsion is unnecessary in this case, which is crucial in applying local
Lefschetz theorem to study the Betti numbers of $\mhd$\ carried out
in [Li3].

\pro{Lemma 3.7}
For any precompact neighborhood $\cc\subb\nqp$, there is a constant
$N$\ such that for any $d\geq N$\ and $(I,d)$-generic
$H\in \cc$, the divisor $\Lambda\sub\mhdb$\
consisting of non-locally free sheaves is Cartier.
\endpro

\proof
Let $k$\ be the number of torsion elements (including 0) in
$H^2(X,\ZZ)$\ and let $\pi\mh\tilx\to X$\ be
a $k$-fold (unramified) covering so that $H^2(\tilx,\ZZ)$\ is torsion free.
Let $\tilde I=\pi\sta I$\ and $\tilde\cc\sub \nqp(\tilx)$\ be the image in
$H^2(\tilx,\QQ)$\ under pullback $H^2(X;\QQ)\to H^2(\tilx;\QQ)$.
$\tilde \cc\subb\nqp(\tilx)$\ is precompact.
Therefore there is an $N$\ such that for any $\tilde d\geq N$\
and $H'\in\tilde\cc$, the moduli scheme
$\overline{\m}_{H\pri}(\tilde I,kd)$\ and the corresponding Quot-scheme
$\qu_{H\pri}(\tilde I,kd)\ss$\ satisfy the conclusion of Lemma 3.1 and 3.2.
Let $\tilde\Lambda\sub\overline{\m}_{H'}(\tilde I,kd)$\
be the closed subset of all non-locally free sheaves. Since
$\pi\mh \tilx\to X$\ is a finite unramified covering, the morphism
$$\overline{\pi}: \mhdb\lra \overline{\m}_{\tilde H}(\tilde I,kd)
$$
induced by sending $\e\in\mhdb$\ to $\pi\sta(\e)$\ is an immersion, where
$\tilde H=\pi\sta H$. In
particular, $\Lambda=\bl\overline{\pi}\br^{-1}\bl\tilde\Lambda\br$.
We claim that $\tilde\Lambda$\ is Cartier near
$\overline{\pi}(\Lambda)$. To establish this, following the proof of
Lemma 3.6, we only need to show that
for any split semistable sheaf $\e=\f_1\oplus\f_2$\ in $\Lambda$,
$\pi\sta(\f_1)\oplus \pi\sta(\f_2)$\ can be deformed within
the set of splitting semistable sheaves in $\overline{\m}_{\tilde
H}(\tilde I, kd)$\ to $\f\oplus\f$. But this is always
possible because $c_1(\f_1)-c_1(\f_2)$\ is torsion implies
$c_1\bl\pi\sta(\f_1)\br=c_1\bl\pi\sta(\f_2)\br$.
Therefore $\tilde\Lambda$\ is Cartier near
$\overline{\pi}(\Lambda)$\ and hence
$\Lambda=\overline{\pi}{-1}(\tilde\Lambda)$\ is Cartier. This
completes the proof of the lemma.
\qed

We now state and prove the main theorem of this note.

\pro{Theorem 3.8}
Suppose $H^2(X,\ZZ)$\ has no torsions. Then
for any precompact $\cc\subb\nqp$, there is an integer $N$\
depending on $(X,I,\cc)$\ such that for any $d\geq N$\ and
$(I,d)$-generic $H\in\cc$, the homomorphism
$$\overline{\vphi}:\picxs\oplus \ZZ\lra \picmdb\tz
$$
constructed in (1.13) has finite kernel and cokernel.
\endpro

\pro{Corollary 3.9}
Suppose $H^2(X,\ZZ)$\ has no torsions. Then
for any precompact $\cc\subb\nqp$, there is an integer $N$\
depending on $(X,I,\cc)$\ such that for any $d\geq N$\ and
$H\in\cc$, the homomorphism
$${\vphi}:\picxs\oplus \ZZ\lra \pic\bl\mhd\br\tz
$$
has finite kernel and cokernel.
\endpro

\proof
By enlarging $\cc\subb\nqp$\ if necessary, we can assume $\cc$\ contains
$\QQ$-polarization $\tilde H$\ that is $(I,d)$-generic for any $d\geq
0$. We then choose $N$\ so that all previous results holds for this $N$.
Now for any $d\geq N$, let $\tilde H\in\cc$\ be $(I,d)$-generic. Then by
Lemma 3.5, $\pic\bl\mhd\br$\ is isomorphic to $\pic\bl\m_{\tilde
H}(I,d)\br$, which is isomorphic to $\pic\bl\overline{\m}_{\tilde
H}(I,d)\br$\ by Lemma 3.6. On the other hand, the homomorphism $\vphi$\
above certainly commutes with
$$\overline{\vphi}: \picxs\oplus \ZZ\lra \pic\bl\overline{\m}_{\tilde
H}(I,d)\br
$$
under these isomorphisms. Therefore, $\vphi$\ has finite kernel and
cokernel because $\overline{\vphi}$\ does. This completes the proof of
the corollary.
\qed

Before we prove the theorem, we need two more technical lemmas. The
first concerns the Hodge decomposition of $H^2\bl\mhd;\RR\br$.
Since we know $H^2\bl\mhd;\RR\br$\ explicitly, we will
determine its Hodge decomposition with the aid of
$$g: W\lra \mhd
$$
that is constructed in the beginning of \S 2 based on a $\mu$-stable
rank two locally free sheaf $\e$\ with $\det \e=I$\ and $c_2(\e)=d-2$.
Let
$$g^{\text{pic}}: \pic\bl\mhd\br\lra\pic\bl W\br
$$
be the induced homomorphism.

\pro{Lemma 3.9}
For any precompact $\cc\subb\nqp$, there is an integer $N$\
depending on $(X,I,\cc)$\ such that for any $d\geq N$, $H\in\cc$\
and $g: W\lra \mhd$\ as before, the image
$$c_1\circ g^{\text{pic}}:\pic\bl\mhd\br\otimes_{\ZZ}\QQ\lra
H^2(W;\QQ)
$$
is isomorphic to
$$H^{1,1}(\xxx;\RR)^{\sigma}\cap H^2(\xxx;\QQ)\oplus \QQ.
$$
\endpro

\proof
Let $\Lambda$\ be the image of $c_1\circ g^{\text{pic}}$.
{}From the gauge theory and the proof of Lemma 3.3 in [Li3],
the image of
$$g\sta: H^2\bl\mhd;\RR\br\lra H^2(W;\RR)
$$
is spanned by the images
$$p_{W\ast}\bl c_2(\g)\cup p_X\sta(\cdot)\br: H^2(X;\RR)\lra
H^2\bl W;\RR\br,
$$
$$p_{W\ast}\bl c_3(\g)\cup p_X\sta(\cdot)\br: H^0(X;\RR)\lra
H^2\bl W;\RR\br
$$
and the image
$$\wedge^2 p_{W\ast}\bl c_2(\g)\cup p_X\sta(\cdot)\br:
\wedge^2 H^1(X;\RR)\lra H^1(W;\RR).
$$
where $\g$\ is the sheaf on $W\times X$\ constructed in \S 2 and $p_W,
p_X$\ are projections of $W\times X$.
Based on exact sequences (2.3) and (2.4), it is
straight forward to check that the direct sum of these homomorphisms
$$\mu: H^2(X;\RR)\oplus \wedge^2 H^1(X;\RR)\oplus H^0(X;\RR)
\lra H^2(W;\RR)
$$
is injective. Since $c_i(\g)$\ is an integral class of pure Hodge type, $\mu$\
preserves the integer lattice as well as the obvious Hodge
structures on both sides. In particular,
$$\Lambda\sub \text{Im}\,\bl\mu\br\cap\bl H^{1,1}(W;\RR)\cap H^2(W;\QQ)\br
\tag 3.2
$$
and the later is isomorphic to the direct sum of
$$H^{0}(X;\QQ)
$$
$$H^{1,1}(X;\RR)\cap H^2(X;\QQ)
$$
and
$$\bl H^{1,0}(X;\RR)\wedge H^{0,1}(X;\RR) \br
\cap\bl \wedge^2 H^1(X;\QQ)\br.
$$
On the other hand, by our construction of the homomorphism
$$\vphi_{\QQ}:\picxs\timezq\oplus \QQ\lra
\pic\bl\mhd\br\timezq
$$
and the formulas of their Chern classes (see (1.6), (1.7) and (1.8)),
we see that the homomorphism
$$H^{1,1}(\xxx;\RR)^{\sigma}\cap H^2(\xxx;\QQ)\oplus\QQ
\lra H^2(W,\QQ)
\tag 3.3
$$
induced by
$$g^{\text{pic}}\circ\vphi_{\QQ}:
\picxs\otimes_{\ZZ}\QQ\oplus\QQ\lra \pic(W)\otimes_{\ZZ}\QQ
$$
is injective. Clearly, the image of (3.3) is contained in $\Lambda$.
Therefore, if we can show that the left hand side of (3.3) is isomorphic
to the right hand side of (3.2), then for dimensional reason the
equality in (3.2) must hold and then the lemma follows immediately.

We now show that
$$\align
H^{1,1}(\xxx;\RR)^{\sigma}\cap H^2(\xxx;\QQ)
&\cong H^{1,1}(X;\RR)\cap H^2(X;\QQ)\\
&\oplus\bl H^{1,0}(X;\RR)\wedge H^{0,1}(X;\RR) \br
\cap\bl \wedge^2 H^1(X;\QQ)\br.
\tag 3.4
\endalign
$$
We will use the Kunneth decomposition
$$H^2(\xxx;\RR)=\operatornamewithlimits{\bigoplus}_{i+j=2}
H^i(X;\RR)\otimes H^j(X;\RR).
$$
It is straightforward to check that an element
$$v=\sum_{i=1}^n \bl b_i(1\otimes\alpha_i)+b_i\pri(\alpha_i\otimes 1)\br
+\sum_{i,j=1}^ma_{ij}(\beta_i\otimes \beta_j)
$$
in $H^2(\xxx;\RR)$\ is $\sigma$\ invariant if and only if
$$b_i=b_i\pri\quad\and\quad a_{ij}=-a_{ji},
$$
where $\{\alpha_i\}^n$\ and $\{\beta_j\}^m$\ are basis of
$H^2(X;\RR)$\ and $H^1(X;\RR)$\ respectively. Hence
$$H^2(\xxx;\RR)^{\sigma}\cong H^2(X;\RR)\oplus
\wedge^2 H^1(X;\RR).
$$
We can derive a similar formula for $H^{1,1}(\xxx;\RR)^{\sigma}$\ and
$H^{1,1}(\xxx;\QQ)^{\sigma}$. Combined, we get the desired identity
(3.4). This completes the proof of Lemma 3.9.
\qed

The next lemma states that the restriction of $c_1(L)$\ to
$H^2\bl\mhd;\ZZ)$, where $L$\ is a line bundle on $\mhdb$,
determines $c_1(L)$\ completely. This is crucial to our
study since we only know the structure of $H^2\bl\mhd;\RR)$\ rather
than $H^2\bl\mhdb;\RR\br$.

\pro{Lemma 3.10}
Let the notation be as in Theorem 3.8. Then there is an $N$\ depending
on $(X,I,\cc)$\ such that for any $d\geq N$\ and $(I,d)$-generic
$H\in\cc$\ the following holds:
Let $L$\ be a line bundle on $\mhdb$\ such that the restriction of
its Chern class to $H^2\bl\mhd;\ZZ\br$\ is trivial, then $L$\ is
derived from a representation
$$\rho: \pi_1\bl\mhdb\br\lra S^1.
$$
\endpro

\proof
Let $Z$\ be a desingularization of $\mhdb$\ with projection
$\pi\mh Z\to\mhdb$\ and let $D_1,\cdots,D_k$\ be irreducible
components of the exceptional divisor of $\pi$, which we assume
is normal crossing. Then the kernel of the composition
$$H^2(Z;\ZZ)\mapright{\pi\sta} H^2\bl\mhdb;\ZZ\br\mapright
{\text{rest.}} H^2\bl\mhdb\lreg;\ZZ\br,
$$
where $\mhdb\lreg$\ is the regular locus of $\mhdb$, is spanned by
all $c_1(\OO(D_i))$. In particular there are integers $n_1,\cdots,n_k$\
such that
$$c_1\bl\pi\sta L(\hbox{$\sum$} n_i D_i)\br=0\in H^2(Z;\ZZ),
$$
since the restriction to $H^2\bl\mhd;\ZZ\br$\ of $c_1(L)$\ is trivial. Hence
$$\tilde L=\pi\sta L\bl\hbox{$\sum$} n_i D_i)
$$
is induced by a representation $\tilde\rho\mh\pi_1(Z)\to S^1$, or induced
by a local system $\tilde\Sigma$.

It remains to show that $\tilde\Sigma$\ descends to a local system
$\Sigma$\ on $\mhdb$, since then the line bundle induced by $\Sigma$\
is isomorphic to $L$\ on $\mhdb\lreg$\ and then by 3.1,
they are isomorphic on $\mhdb$\ as well. This will prove the lemma.

Instead of studying the descent of $\tilde\Sigma$\ via projection
$\pi\mh Z\to\mhdb$, which requires the knowledge of the singularity of
$\mhdb$, we will work with the quotient $p\mh\qu_H(I,d)\ss\to
\mhdb$\ since virtually all local information of $\mhdb$\ is contained in
this construction. First, since $\pi\mh Z\to\mhdb$\ is an isomorphism
over $\mhdb\lreg$, we obtain a local system $\Sigma\pri$\ on
$\mhdb\lreg$\ and its pull-back $p\sta(\Sigma\pri)$\ on
$p^{-1}\bl\mhdb\lreg\br\sub\qu_H(I,d)\ss$. Clearly, $p\sta(\Sigma\pri)$\
is $G$-equivalent. Let
$$R=\qu_H(I,d)\ss-p^{-1}\bl \mhdb\lreg\br.
$$
By Lemma 3.1, $R$\ has codimension at least 5. Then for any closed $x\in
R$\ and normal slice $S$\ of $R\sub\qu_H(I,d)\ss$\ at $x$,
$\pi_1(S-x)=\{1\}$, because
$\qu_H(I,d)\ss$\ is a local complete intersection scheme
((2) of Lemma 3.2) and $\codim(R)\geq 5$, by the local Lefschetz theorem
on page 155 of [GM].
Therefore $p\sta(\Sigma\pri)$\ extends to a local system, say
$\Sigma''$, that is $G$-equivalent as well, since $\qu_H(I,d)\ss$\
is normal. To show that $\Sigma''$\ descends to $\Sigma$\ on $\mhdb$,
it suffices to show that for any $w\in\qu_H(I,d)\ss$\ with closed
orbit $G\cdot w$, $\stab_w$\ acts trivially on the fiber of $\Sigma''$\
over $w$. But this can be proved using the fact that $H^2(X;\ZZ)$\
has no torsion as we did in the proof of Lemma 3.6. This shows that $\Sigma''$\
descends to a local system $\Sigma$. This completes the proof
of the Lemma.
\qed

\noindent
{\it Proof of Theorem 3.8}.
We let $N$\ be the integer that makes all previous results valid. Let
$d\geq N$\ and let $H\in\cc$\ be $(I,d)$-generic. By Lemma 3.8 for any
line bundle $L$ on $\mhdb$, the restriction of its Chern class
$c_1(L)$\ to $H^2\bl\mhd;\QQ)$\ is contained in the image of
$$c_1\circ\vphi: \picxs\timezq\oplus\QQ\lra H^2\bl\mhd;\QQ\br.
$$
Hence we can find an element $\alpha\in\picxs\oplus\ZZ$\ and an
integer $k$\ such that the restriction of $c_1\bl L^{\otimes k}\otimes
L_{\alpha}\br$\ to $H^2\bl \mhd;\ZZ\br$\ is trivial, where $L_{\alpha}$\
is the line bundle $\overline{\vphi}(\alpha)$\ on $\mhdb$.
Here $k$\ can be chosen to depend only on the induced homomorphism
$$H^2(\xxx;\ZZ)^{\sigma}\oplus\ZZ\lra H^2\bl\mhd;\ZZ[{1\over 12}]\br.
$$
By Lemma 3.10, $L^{\otimes k}\otimes
L_{\alpha}$\ is induced by a representation
$\rho:\pi_1\bl\mhdb\br\lra S^1$, or an element in
$H^1\bl\mhdb;\RR\br$. However, $\mhdb$\ is normal and $\mhdb-\mhd$\
has codimension at least 2, hence
$$H^1\bl\mhdb;\RR\br\lra H^1\bl\mhd;\RR\br
$$
is injective. On the other hand, we know
$$\dim H^1\bl\mhd;\RR\br=h^1(X)
$$
and the restriction of $\vphi$
$$\vphi_0: \pic^0\bl\xxx\br^{\sigma}\lra \pic^0\bl\mhdb\br
$$
has finite kernel since
$$g^{\text{pic}}\circ\vphi_0:\picxs\lra \pic(W)
$$
has finite kernel, by direct check. This implies $\text{Im}(\vphi_0)$\
has dimension $h^1(X)$. Therefore, $\vphi_0$\ is surjective. This
proves that
$$\vphi:\picxs\oplus \ZZ\lra \pic\bl\mhdb\br\tz
$$
has finite cokernel.
\qed

\pro{Remark}
1. The author inclined to believe the condition that $H^2(X;\ZZ)$\ is
torsion free in this theorem is unnecessary, but is unable to remove it
using the current technique.

2. It is interesting to to find a bound $N$\ in this theorem that can be
determined effectively. After the effective bound of O'Grady [OG], we
can find an effective bound that works for all results in this note that
are independent of Lemma 3.3 whose proof relies on the result of C.
Taubes that is not effective. It is interesting to get a purely algebraic
proof of Lemma 3.3.
\endpro

\vskip0.5in
\noindent
Mathematics Department, University of California at Los Angeles, CA 90024,
USA.

\noindent
email: jli\@math.ucla.edu

\bye